\documentclass[a4paper,11pt,twoside]{amsart}

\usepackage[latin1]{inputenc}
\usepackage[T1]{fontenc}
\usepackage{amsmath}
\usepackage{amsthm}
\usepackage{amssymb}
\usepackage[all]{xy}

\newtheorem{thm}{Theorem}[section]
\newtheorem{prop}[thm]{Proposition}
\newtheorem{lem}[thm]{Lemma}
\newtheorem{cor}[thm]{Corollary}

\numberwithin{equation}{section}

\theoremstyle{definition}
\newtheorem{definition}[thm]{Definition}
\newtheorem{remark}[thm]{Remark}
\newtheorem{ex}[thm]{Example}

\newcommand{\supp}{\operatorname{Supp}}

\newcommand{\im}{\operatorname{im}}

\newcommand{\Db}{{\rm D}^{\rm b}}
\newcommand{\D}{{\rm D}}

\newcommand{\Br}{{\rm Br}}

\newcommand{\Pic}{{\rm Pic}}

\newcommand{\coh}{{\cat{Coh}}}
\newcommand{\qcoh}{{\cat{QCoh}}}

\newcommand{\Hom}{{\rm Hom}}

\newcommand{\ob}{{\rm Ob}}

\newcommand{\iso}{\cong}

\newcommand{\id}{{\rm id}}

\newcommand{\comp}{\circ}
\newcommand{\rest}[1]{|_{#1}}
\newcommand{\dual}{^{\vee}}
\newcommand{\mono}{\hookrightarrow}
\newcommand{\epi}{\twoheadrightarrow}
\newcommand{\mor}[1][]{\xrightarrow{#1}}
\newcommand{\isomor}{\mor[\sim]}
\newcommand{\abs}[1]{\lvert#1\rvert}
\newcommand{\cp}[1]{#1^{\bullet}}
\newcommand{\cat}[1]{\begin{bf}#1\end{bf}}

\newcommand{\FM}[1]{\Phi_{#1}}
\newcommand{\aFM}[1]{\Psi_{#1}}
\newcommand{\R}{\mathbf{R}}
\newcommand{\lotimes}{\stackrel{\mathbf{L}}{\otimes}}
\newcommand{\lboxtimes}{\stackrel{\mathbf{L}}{\boxtimes}}

\newcommand{\K}{{\rm K}}

\newcommand{\cal}{\mathcal}
\newcommand{\ka}{{\cal A}}
\newcommand{\kb}{{\cal B}}

\newcommand{\ke}{{\cal E}}
\newcommand{\kf}{{\cal F}}
\newcommand{\kg}{{\cal G}}

\newcommand{\kk}{{\cal K}}

\newcommand{\ko}{{\cal O}}

\newcommand{\kx}{{\cal X}}
\newcommand{\ky}{{\cal Y}}

\newcommand{\NN}{\mathbb{N}}
\newcommand{\ZZ}{\mathbb{Z}}

\newcommand{\PP}{\mathbb{P}}
\newcommand{\sHom}{\cal{H}om}
\newcommand{\ds}{\omega^{\circ}}

\begin{document}

\title[Twisted Fourier-Mukai functors]{Twisted Fourier-Mukai functors}

\author{Alberto Canonaco and Paolo Stellari}

\address{A.C.: Dipartimento di Matematica ``F. Casorati'', Universit{\`a}
degli Studi di Pavia, Via Ferrata 1, 27100 Pavia, Italy}
\email{alberto.canonaco@unipv.it}

\address{P.S.: Dipartimento di Matematica ``F.
Enriques'', Universit{\`a} degli Studi di Milano, Via Cesare Saldini
50, 20133 Milano, Italy}
\email{stellari@mat.unimi.it}



\begin{abstract} Due to a theorem by Orlov every exact fully
faithful functor between the bounded derived categories of
coherent sheaves on smooth projective varieties is of
Fourier-Mukai type. We extend this result to the case of bounded
derived categories of twisted coherent sheaves and at the same
time we weaken the hypotheses on the functor. As an application we
get a complete description of the exact functors between the
abelian categories of twisted coherent sheaves on smooth
projective varieties.\end{abstract}

\maketitle

\section{Introduction}\label{sec:introduction}

If $X$ and $Y$ are smooth projective varieties, an exact functor
$F:\Db(X)\to\Db(Y)$ between the corresponding bounded derived
categories of coherent sheaves is of \emph{Fourier-Mukai type} if
there exists $\ke\in\Db(X\times Y)$ and an isomorphism of functors
$F\iso\FM{\ke}$, where, denoting by $p:X\times Y\to Y$ and
$q:X\times Y\to X$ the natural projections,
$\FM{\ke}:\Db(X)\to\Db(Y)$ is the exact functor defined by
\begin{eqnarray}\label{eqn:FMT}
\FM{\ke}:=\R p_*(\ke\lotimes q^*(-)).
\end{eqnarray}
Such a complex $\ke$ is called a \emph{kernel} of $F$.

The importance of functors of this type in geometric contexts
cannot be overestimated. Indeed, all meaningful geometric functors
are of Fourier-Mukai type and conjecturally the same is true for
every exact functor from $\Db(X)$ to $\Db(Y)$. As a first evidence
for the truth of this conjecture, in the fundamental paper
\cite{Or1}, Orlov proved that any exact fully faithful functor
from $\Db(X)$ to $\Db(Y)$ which admits a left adjoint is of
Fourier-Mukai type. Moreover its kernel is uniquely determined up
to isomorphism.

Since the publication of \cite{Or1}, some significant improvements
were obtained. The main one is due to Kawamata (\cite{Ka}), who
extended this result to the case of smooth quotient stacks. His
proof partially follows Orlov's original one but at some crucial
points new deep ideas are needed. It is also worth noticing that,
due to the results in \cite{BB}, every exact functor
$F:\Db(X)\to\Db(Y)$ admits a left adjoint (see Remark \ref{rmk:adj}
below).

In recent years some attention was paid to the case of \emph{twisted
varieties} (i.e.\ pairs $(X,\alpha)$, where $X$ is a smooth
projective variety and $\alpha$ is an element in the Brauer group of
$X$). Since \cite{C} appeared, it has been proved that some results
from the untwisted setting can be generalized to the case of twisted
derived categories. For example, if $M$ is a K3 surface and a moduli
space of stable sheaves on a K3 surface $X$, then there exist
$\alpha$ in the Brauer group of $M$ and an equivalence between
$\Db(X)$ and $\Db(M,\alpha)$, the bounded derived category of
$\alpha$-twisted coherent sheaves on $M$. This was first proved by
C\u{a}ld\u{a}raru (\cite{C,C2}) and then generalized in
\cite{L,L1,Y} and \cite{HS1,HS2}. Nevertheless a question remained
open:

\medskip

\noindent\centerline{\it Are all equivalences between the bounded
derived categories of twisted} \centerline{\it coherent sheaves on
smooth projective varieties of Fourier-Mukai type?}

\medskip

\noindent As before, given two twisted varieties $(X,\alpha)$ and
$(Y,\beta)$, a functor $F:\Db(X,\alpha)\to\Db(Y,\beta)$ is of
Fourier-Mukai type if there exist $\ke\in\Db(X\times
Y,\alpha^{-1}\boxtimes\beta)$ and an isomorphism of functors
$F\iso\FM{\ke}$, where $\FM{\ke}$ is again defined as in
\eqref{eqn:FMT}.

A complete answer to the previous question comes as an easy
corollary of the following theorem which is the main result of this
paper:

\begin{thm}\label{thm:main} Let $(X,\alpha)$ and $(Y,\beta)$ be
twisted varieties and let $F:\Db(X,\alpha)\to\Db(Y,\beta)$ be an
exact functor such that, for any $\kf,\kg\in\coh(X,\alpha)$,
\begin{eqnarray}\label{eqn:hyp}
\Hom_{\Db(Y,\beta)}(F(\kf),F(\kg)[j])=0\;\;\mbox{if }j<0.
\end{eqnarray}
Then there exist $\ke\in\Db(X\times Y,\alpha^{-1}\boxtimes\beta)$
and an isomorphism of functors $F\iso\FM{\ke}$. Moreover, $\ke$ is
uniquely determined up to isomorphism.\end{thm}

Few comments about the relevance of the previous result are in order
here. First of all observe that any full functor satisfies
\eqref{eqn:hyp}. This means that Theorem \ref{thm:main} gives a
substantial improvement of Orlov's result. As a consequence, we will
observe that also the hypotheses in Kawamata's result (\cite{Ka})
can be weakened (see Remark \ref{Rmk:kaw}).

Our proof of Theorem \ref{thm:main} was inspired by \cite{Ka} and
\cite{Or1} although different approaches are needed in many
crucial points. In particular, the idea to use extensively
convolutions of bounded complexes comes from \cite{Or1}.

In Section \ref{sec:exfun} we apply Theorem \ref{thm:main} to
describe exact functors between the abelian categories of twisted
coherent sheaves. In particular we deduce a Gabriel-type result for
twisted varieties.

\medskip

\noindent{\bf Notations.} We will work over a fixed field $\K$. All
triangulated and abelian categories and all exact functors will be
assumed to be $\K$-linear. For an abelian category $\cat{A}$ we will
denote by $\D(\cat{A})$ the derived category of $\cat{A}$. An object
$\cp{C}$ of $\D(\cat{A})$ is a complex in $\cat{A}$, i.e.\ it is
given by a collection of objects $C^i$ and morphisms $d^i:C^i\to
C^{i+1}$ of $\cat{A}$ such that $d^{i+1}\comp d^i=0$. The bounded
derived category of $\cat{A}$ is the full subcategory $\Db(\cat{A})$
of $\D(\cat{A})$ with objects the complexes $\cp{C}$ such that
$C^i=0$ for $\abs{i}\gg0$. If there is no ambiguity, we will usually
write $C$ instead of $\cp{C}$. If $\cat{B}$ is another abelian
category, every exact functor $G:\cat{A}\to\cat{B}$ trivially
induces exact functors of triangulated categories
$\D(G):\D(\cat{A})\to\D(\cat{B})$ and
$\Db(G):\Db(\cat{A})\to\Db(\cat{B})$. Recall that an abelian
category $\cat{A}$ is of finite homological dimension if there
exists an integer $l$ such that, for any $i>l$ and any
$A,B\in\ob(\cat{A})$, $\Hom_{\Db(\cat{A})}(A,B[i])=0$; if $N\in\NN$
is the least such integer $l$, then $\cat{A}$ is said to be of
homological dimension $N$.

\section{Boundedness and ample sequences}\label{subsec:amseq}

For a smooth projective variety $X$ consider the cohomology group
$H^2_{\mathrm{\acute{e}t}}(X,\ko_X^*)$ in the \'{e}tale topology.
Any $\alpha\in H^2_{\mathrm{\acute{e}t}}(X,\ko_X^*)$ can be
represented by a \v{C}ech 2-cocycle on an \'{e}tale cover
$\{U_i\}_{i\in I}$ of $X$ using sections
$\alpha_{ijk}\in\Gamma(U_i\cap U_j\cap U_k,\ko^*_X)$. An
\emph{$\alpha$-twisted quasi-coherent sheaf} $\kf$ consists of a
pair $(\{\kf_i\}_{i\in I},\{\varphi_{ij}\}_{i,j\in I})$, where
$\kf_i$ is a quasi-coherent sheaf on $U_i$ and
$\varphi_{ij}:\kf_j|_{U_i\cap U_j}\to\kf_i|_{U_i\cap U_j}$ is an
isomorphism such that $\varphi_{ii}=\mathrm{id}$,
$\varphi_{ji}=\varphi_{ij}^{-1}$ and
$\varphi_{ij}\circ\varphi_{jk}\circ\varphi_{ki}=\alpha_{ijk}\cdot\mathrm{id}$.

The category of $\alpha$-twisted quasi-coherent sheaves on $X$
will be denoted by $\qcoh(X,\alpha)$. An $\alpha$-twisted
quasi-coherent sheaf $(\{\kf_i\}_{i\in I},\{\varphi_{ij}\}_{i,j\in
I})$ is an \emph{$\alpha$-twisted coherent sheaf} if $\kf_i$ is
coherent for any $i\in I$. We write $\coh(X,\alpha)$ for the
abelian category of $\alpha$-twisted coherent sheaves and
$\Db(X,\alpha):=\Db(\coh(X,\alpha))$ for the bounded derived
category of $\coh(X,\alpha)$. The \emph{Brauer group} of $X$ is
the group $\Br(X)$ consisting of all $\alpha\in
H^2_{\mathrm{\acute{e}t}}(X,\ko_X^*)$ such that $\coh(X,\alpha)$
contains a locally free $\alpha$-twisted coherent sheaf (actually,
due to \cite{dJ}, $\Br(X)$ coincides with
$H^2_{\mathrm{\acute{e}t}}(X,\ko_X^*)$).

Let $X$ and $Y$ be smooth projective varieties and let $f:X\to Y$ be
a morphism. The following derived functors are defined:
$-\lotimes-:\Db(X,\alpha)\times\Db(X,\alpha')\to\Db(X,\alpha\cdot\alpha')$,
$\R f_*:\Db(X,f^*(\beta))\to\Db(Y,\beta)$ and
$\mathbf{L}f^*:\Db(Y,\beta)\to\Db(X,f^*(\beta))$, where
$\alpha,\alpha'\in\Br(X)$ and $\beta\in\Br(Y)$ (see \cite[Thm.\
2.2.4, Thm.\ 2.2.6]{C}). For the rest of this paper $(X,\alpha)$ and
$(Y,\beta)$ will denote two twisted varieties as in Theorem
\ref{thm:main}.

\begin{remark}\label{rmk:adj} (i) If $(X,\alpha)$ is a twisted variety and $X$
has dimension $n$, then $\coh(X,\alpha)$ has homological dimension $n$. To
prove this claim, one can proceed as in the untwisted case (see \cite[Prop.\
3.12]{Hu1}), using the fact that the functor $S(-)=(-)\otimes\omega_X[n]$ is
the Serre functor of $\Db(X,\alpha)$.

(ii) If $(X,\alpha)$ and $(Y,\beta)$ are twisted varieties, then any
exact functor $G:\Db(X,\alpha)\to\Db(Y,\beta)$ has a left adjoint
$G^*:\Db(Y,\beta)\to\Db(X,\alpha)$. Indeed, it is proved in \cite{P}
(generalizing ideas from \cite{BB} and \cite{R}) that any
cohomological functor of finite type is representable. Hence, for
any $\kf\in\Db(Y,\beta)$ the functor $\Hom_{\Db(Y,\beta)}(G(-),\kf)$
is representable by a unique $\ke\in\Db(X,\alpha)$. Setting
$G'(\kf):=\ke$, by the Yoneda Lemma we get a functor which is right
adjoint to $G$. Since $\Db(X,\alpha)$ and $\Db(Y,\beta)$ have Serre
functors, $G$ has also a left adjoint $G^*$.\end{remark}

\begin{definition}\label{def:ample} Given an abelian category
$\cat{A}$ with finite dimensional $\Hom$'s, a subset
$\{P_i\}_{i\in\ZZ}\subset\ob(\cat{A})$ is an \emph{ample sequence}
if, for any $B\in\ob(\cat{A})$, there exists an integer $i(B)$ such
that, for any $i\leq i(B)$,
\begin{enumerate}
\item\label{ample1} the natural morphism $\Hom_{\cat{A}}(P_i,B)
\otimes P_i\to B$ is surjective;

\item\label{ample2} if $j\ne0$ then
$\Hom_{\Db(\cat{A})}(P_i,B[j])=0$;

\item\label{ample3} $\Hom_{\cat{A}}(B,P_i)=0$.
\end{enumerate}\end{definition}

\begin{lem}\label{lem:ample} Let $E\in\coh(X,\alpha)$ be a locally free sheaf.
If $\{\ka_k\}_{k\in\ZZ}$ is an ample sequence in $\coh(X)$, then
$\{E\otimes\ka_k\}_{k\in\ZZ}$ is an ample sequence in $\coh(X,\alpha)$. In
particular, if $L\in\coh(X)$ is an ample line bundle, then $\{E\otimes
L^{\otimes k}\}_{k\in\ZZ}$ is an ample sequence.\end{lem}

\begin{proof} Observe that since $\{\ka_k\}_{k\in\ZZ}$ is an ample
sequence, for any $\ke\in\coh(X,\alpha)$ and for $i\ll 0$, there
exists a surjective map
$\Hom_{\coh(X)}(\ka_i,E^\vee\otimes\ke)\otimes\ka_i\epi
E^\vee\otimes\ke$. Then \eqref{ample1} in the previous definition
follows from the fact that the diagram
\[
\xymatrix{\Hom_{\coh(X,\alpha)}(E\otimes \ka_i,\ke)
\otimes E\otimes \ka_i \ar[r] \ar[d]_{\iso} & \ke \\
\Hom_{\coh(X)}(\ka_i,E\dual\otimes\ke)\otimes E\otimes \ka_i
\ar@{->>}[r] & E\otimes E\dual\otimes\ke \ar@{->>}[u]}
\]
commutes. Analogously,
$\Hom_{\Db(X,\alpha)}(E\otimes\ka_i,\ke[j])\cong
\Hom_{\Db(X)}(\ka_i,E^\vee\otimes\ke[j])=0$ and
$$\Hom_{\coh(X,\alpha)}(\ke,E\otimes\ka_i)\cong\Hom_{\coh(X)}(\ke\otimes
E^\vee,\ka_i)=0,$$ for $i\ll0$ and $j\neq 0$. This proves that
\eqref{ample2} and \eqref{ample3} hold true. The second part of the
lemma follows from the easy fact that $\{L^{\otimes k}\}_{k\in\ZZ}$
is an ample sequence in $\coh(X)$.\end{proof}

Recall that, given two abelian categories $\cat{A}$ and $\cat{B}$,
a functor $G:\Db(\cat{A})\mor\Db(\cat{B})$ is \emph{bounded} if
there exist $a\in\ZZ$ and $n\in\NN$ such that $H^i(G(A))=0$ for
any $A\in\ob(\cat{A})$ and any $i\not\in[a,a+n]$.

\begin{prop}\label{prop:bound} Let $(X,\alpha)$ and $(Y,\beta)$ be
twisted varieties and assume that $G:\Db(X,\alpha)\mor\Db(Y,\beta)$
is an exact functor. Then $G$ is bounded.\end{prop}

\begin{proof} Due to Lemma \ref{lem:ample}, given a
locally free sheaf $E\in\coh(Y,\beta)$ and a very ample line bundle
$L\in\coh(Y)$ (defining an embedding $Y\mono\PP^N$), the set
$\{E\otimes L^{\otimes k}\}_{k\in\ZZ}$ is an ample sequence in
$\coh(Y,\beta)$. For $k<0$, Beilinson's resolution (\cite{B}),
pulled back to $Y$, yields an isomorphism in $\Db(Y)$
\begin{eqnarray}\label{eqn:Beilinson}
L^{\otimes k}\iso\{V^k_N\otimes\ko_Y\to V^k_{N-1}\otimes
L\to\ldots\to V^k_0\otimes L^{\otimes N}\},
\end{eqnarray}
where $V^k_i:=H^N(\PP^N,\Omega^i_{\PP^N}(i+k-N))$. In particular, $E\otimes
L^{\otimes k}\iso\cp{C_k}$ in $\Db(Y,\beta)$ where $C_k^i=0$, for $\abs{i}>N$,
and each $C_k^i$ is a finite direct sum of terms of the form $E\otimes L^j$ for
$0\le j\le N$. This implies that $\{G^*(E\otimes L^{\otimes k})\}_{k<0}$ is
bounded in $\Db(X,\alpha)$, where $G^*$ is the left adjoint of $G$ (see Remark
\ref{rmk:adj}(ii)).

This is enough to conclude that $G$ is bounded. Indeed we can reason in the
following rather standard way. Given $\ka\in\coh(X,\alpha)$ and $i\in\ZZ$, it
is easy to see that $H^i(G(\ka))=0$ is implied by $\Hom_{\Db(Y,\beta)}(E\otimes
L^{\otimes k},G(\ka)[i])=0$ for $k\ll 0$. Choosing $m$ such that
$H^j(G^*(E\otimes L^{\otimes k}))=0$ for $\abs{j}\ge m$ and for $k<0$ and
denoting by $n$ the homological dimension of $\coh(X,\alpha)$ (see Remark
\ref{rmk:adj}(i)), it is clear that
$$\Hom_{\Db(Y,\beta)}(E\otimes L^{\otimes k},G(\ka)[i])\cong
\Hom_{\Db(X,\alpha)}(G^*(E\otimes L^{\otimes k}),\ka[i])=0$$ for
$\abs{i}>n+m$ and for $k<0$.\end{proof}

The following easy lemma will be used in the forthcoming sections.

\begin{lem}\label{lem:serrerel} Let $(X,\alpha)$ and $(Y,\beta)$ be
twisted varieties, let $\ke\in\Db(X\times Y,\alpha^{-1}\boxtimes\beta)$ and let
$l\in\ZZ$. If $E\in\coh(X,\alpha)$ is locally free and $L\in\coh(X)$ is ample,
then $H^l(\ke)=0$ if $H^l(\FM{\ke}(E\otimes L^{\otimes k}))=0$ for any $k\gg
0$.\end{lem}

\begin{proof} Fix a locally free sheaf $F\in\coh(Y,\beta)$ and define
$\ka:=\ke\otimes p^*(F\dual)\otimes q^*(E)$, where $p:X\times Y\to Y$ and
$q:X\times Y\to X$ are the natural projections. If $H^j(\ka)\ne0$, $\R^i
p_*(H^j(\ka)\otimes q^*(L^{\otimes k}))=0$ for any $k\gg 0$ if and only if
$i\neq 0$ (for a proof of this well-known fact see, for example, \cite{Ha},
Chapter III, Theorem 8.8). Hence, using the spectral sequence
\[
\R^i p_*(H^j(\ka)\otimes q^*(L^{\otimes k}))\Longrightarrow
H^{i+j}(\FM{\ka}(L^{\otimes k}))
\]
we deduce that $H^l(\ka)=0$ if $H^l(\FM{\ka}(L^{\otimes k}))=0$, for
any $k\gg 0$.

It is obvious that $H^l(\ka)=0$ if and only if $H^l(\ke)=0$. Hence
the result is proved once we show that $H^l(\FM{\ka}(L^{\otimes
k}))=0$ if and only if $H^l(\FM{\ke}(E\otimes L^{\otimes k}))=0$. By
the Projection Formula
\[\FM{\ka}(L^{\otimes k})\iso\R p_*(\ka\otimes
q^*(L^{\otimes k})) \iso\R p_*(\ke\otimes p^*(F\dual)\otimes
q^*(E)\otimes q^*(L^{\otimes k}))\iso\FM{\ke}(E\otimes L^{\otimes
k})\otimes F\dual
\]
which yields the desired conclusion.\end{proof}

\section{Convolutions and isomorphisms of functors}

In this section we recall few results about convolutions of bounded
complexes and we use them to study the existence of isomorphisms of
exact functors.

\subsection{Convolutions} Recall that a bounded complex in a
triangulated category $\cat{D}$ is a sequence of objects and
morphisms in $\cat{D}$
\begin{equation}\label{eqn:complex}
A_{m}\mor[d_{m}]A_{m-1}\mor[d_{m-1}]\cdots\mor[d_{1}]A_0
\end{equation}
such that $d_{j}\circ d_{j+1}=0$ for $0<j<m$. Following the
terminology of \cite{Ka}, a \emph{right convolution} of
\eqref{eqn:complex} is an object $A$ together with a morphism
$d_0:A_0\to A$ such that there exists a diagram in $\cat{D}$
\[\xymatrix{A_m \ar[rr]^{d_m} \ar[dr]_{\id} \ar@{}[drr]|{\circlearrowright}
& & A_{m-1} \ar[rr]^{d_{m-1}} \ar[dr] \ar@{}[drr]|{\circlearrowright} & &
\cdots \ar[rr]^{d_2} & & A_1 \ar[rr]^{d_1} \ar[dr]
\ar@{}[drr]|{\circlearrowright} & & A_0 \ar[dr]_{d_0} \\
 & A_m \ar[ru] & & C_{m-1} \ar[ll]^{[1]} \ar[ru] & & \cdots
\ar[ll]^{[1]} & & C_1 \ar[ll]^{[1]} \ar[ru] & & A, \ar[ll]^{[1]}}\]
where the triangles marked with a $\circlearrowright$ are
commutative and the other triangles are distinguished (such an
object $A$ is called instead a left convolution of
\eqref{eqn:complex} in \cite{Or1}). In a completely dual way, a
\emph{left convolution} of \eqref{eqn:complex} is an object $A'$
together with a morphism $d_{m+1}:A'\to A_m$ such that there exists
a diagram in $\cat{D}$
\[\xymatrix{& A_m \ar[rr]^{d_m} \ar[dr] \ar@{}[drr]|{\circlearrowright}
& & A_{m-1} \ar[rr]^{d_{m-1}} \ar[dr] & & \cdots \ar[rr]^{d_2}
\ar[dr] \ar@{}[drr]|{\circlearrowright} & & A_1
\ar@{}[drr]|{\circlearrowright} \ar[rr]^{d_1} \ar[dr]
& & A_0 \\
A'\ar[ru]_{d_{m+1}} & & C'_{m-1} \ar[ll]^{[1]} \ar[ru] & & \cdots
\ar[ll]^{[1]} & & C'_1 \ar[ll]^{[1]} \ar[ru] & & A_0. \ar[ll]^{[1]}
\ar[ru]_{\id} &}\]

\begin{remark}\label{rmk:convfun}
Assume that $d_0:A_0\to A$ (respectively $d_{m+1}:A'\to A_m$) is a
right (respectively left) convolution of \eqref{eqn:complex}. If
$\cat{D}'$ is another triangulated category and
$G:\cat{D}\to\cat{D}'$ is an exact functor, then it is obvious from
the definitions that $G(d_0):G(A_0)\to G(A)$ (respectively
$G(d_{m+1}):G(A')\to G(A_m)$) is a right (respectively left)
convolution of
\[G(A_{m})\mor[G(d_{m})]G(A_{m-1})\mor[G(d_{m-1})]\cdots\mor[G(d_{1})]G(A_0).\]
\end{remark}

In general a (right or left) convolution of a complex need not
exist, and it need not be unique up to isomorphism when it exists,
but we have the following two results, which will be constantly used
in the rest of this paper:

\begin{lem}\label{lem:conv1} {\bf (\cite{Ka}, Lemmas 2.1 and 2.4.)}
Let \eqref{eqn:complex} be a complex in $\cat{D}$ satisfying
\begin{equation}\label{eqn:condconv}
\Hom_{\cat{D}}(A_a,A_b[r])=0\mbox{ for any $a>b$ and $r<0$}.
\end{equation}
Then \eqref{eqn:complex} has right and left convolutions and they
are uniquely determined up to isomorphism (in general non
canonical).
\end{lem}

\begin{lem}\label{lem:conv2}
Let
\[\xymatrix{A_m\ar[r]^{d_m}\ar[d]^{f_m}&A_{m-1}\ar[r]^{d_{m-1}}\ar[d]^{f_{m-1}}
&\cdots\ar[r]^{d_2}&A_1\ar[r]^{d_1}\ar[d]^{f_1}&A_0\ar[d]^{f_0}\\
B_m\ar[r]^{e_m}&B_{m-1}\ar[r]^{e_{m-1}}&\cdots\ar[r]^{e_2}&B_1\ar[r]^{e_1}&B_0}
\]
be a morphism of complexes both satisfying \eqref{eqn:condconv} and
such that
\[\Hom_{\cat{D}}(A_a,B_b[r])=0\mbox{ for any $a>b$ and $r<0$}.\] Assume that
the corresponding right (respectively left) convolutions are of the
form $(d_0,0):A_0\to A\oplus\bar{A}$ and $(e_0,0):B_0\to
B\oplus\bar{B}$ (respectively $(d_{m+1},0):A'\oplus\bar{A}'\to A_m$
and $(e_{m+1},0):B'\oplus\bar{B}'\to B_m$) and that
$\Hom_{\cat{D}}(A_p,B[r])=0$ (respectively
$\Hom_{\cat{D}}(A',B_p[r])=0$) for $r<0$ and any $p$. Then there
exists a unique morphism $f:A\to B$ (respectively $f':A'\to B'$)
such that $f\comp d_0=e_0\comp f_0$ (respectively $e_{m+1}\comp
f'=f_m\comp d_{m+1}$). If moreover each $f_i$ is an isomorphism,
then $f$ (respectively $f'$) is an isomorphism as well.
\end{lem}

\begin{proof}
The first part is a particular case of Lemma 2.3 (respectively Lemma
2.6) of \cite{Ka}. From this it is then straightforward  to deduce
that $f$ (respectively $f'$) is an isomorphism if each $f_i$ is an
isomorphism.
\end{proof}

\begin{ex}\label{ex:convcp}
Let $\cat{D}:=\Db(\cat{A})$ for some abelian category $\cat{A}$ and
let $Z$ be a complex as in \eqref{eqn:complex} and such that every
$A_i$ is an object of $\cat{A}$. Then it is easy to see that a right
(respectively left) convolution of $Z$ (which is unique up to
isomorphism by Lemma \ref{lem:conv1}) is given by the natural
morphism $A_0\to\cp{Z}$ (respectively $\cp{Z}[-m]\to A_m$), where
$\cp{Z}$ is the object of $\Db(\cat{A})$ naturally associated to $Z$
(namely, $Z^i:=A_{-i}$ for $-m\le i\le0$ and otherwise $Z^i:=0$,
with differential $d_{-i}:Z^i\to Z^{i+1}$ for $-m\le i<0$).
\end{ex}

\subsection{Extending isomorphisms of functors}\label{subsec:isom} Let
$\cat{A}$ be an abelian category with finite dimensional $\Hom$'s and assume
that $\{P_i\}_{i\in\ZZ}\subset\ob(\cat{A})$ is an ample sequence.

\begin{lem}\label{lem:resol} Any $A\in\cat{A}$ admits a resolution
\begin{eqnarray}
\cdots\mor A_i^{\oplus k_i}\mor[d_i]A_{i-1}^{\oplus
k_{i-1}}\mor[d_{i-1}]\cdots\mor[d_1]A_0^{\oplus k_0}\mor[d_0] A\mor0,
\end{eqnarray}
where $A_j\in\{P_i\}_{i\in\ZZ}$ and $k_j\in\NN$, for any
$j\in\NN$.\end{lem}

\begin{proof} To prove that such a(n infinite) resolution exists, it
is clearly enough to show that for any $B\in\ob(\cat{A})$ there
exists $P\in\{P_i\}_{i\in\ZZ}$ and a surjective map $P^{\oplus
k}\twoheadrightarrow B,$ for some $k\in\NN$. This follows from
condition (1) in Definition \ref{def:ample}.\end{proof}

\begin{remark}\label{rmk:sp} Consider a resolution of $A\in\ob(\cat{A})$ as in
Lemma \ref{lem:resol} and assume that $\cat{A}$ has finite
homological dimension $N$. Take $m>N$ and consider the bounded
complex
\[
S_m:=\{A_m^{\oplus k_{m}}\mor[d_m]A_{m-1}^{\oplus
k_{m-1}}\mor[d_{m-1}]\cdots\mor[d_1]A_0^{\oplus k_{0}}\}.
\]
If $K_m:=\ker(d_m)$, we have a distinguished triangle in
$\Db(\cat{A})$
\[
K_m[m]\mor[]\cp{S_m}\mor[]A\mor[]K_m[m+1].
\]
Due to the choice of $m$,
$\Hom_{\Db(\cat{A})}(A,K_m[m+1])=\Hom_{\Db(\cat{A})}(A_0,K_m[m])=0$.
Hence $\cp{S_m}\iso A\oplus K_m[m]$ and $S_m$ has a (unique up to
isomorphism) convolution $(d_0,0):A_0^{\oplus k_0}\mor A\oplus
K_m[m]$ (see Example \ref{ex:convcp}).\end{remark}

The following result, whose proof relies on an extensive use of
convolutions, improves \cite[Lemma 6.5]{Ka} and \cite[Prop.\
2.16]{Or1}.

\begin{prop}\label{prop:extending} Let $\cat{D}$ be a
triangulated category and let $\cat{A}$ be an abelian category with finite
dimensional $\Hom$'s and of finite homological dimension. Assume that
$\{P_i\}_{i\in\ZZ}\subseteq\ob(\cat{A})$ is an ample sequence and denote by
$\cat{C}$ the full subcategory of $\Db(\cat{A})$ such that
$\ob(\cat{C})=\{P_i\}_{i\in\ZZ}$. Let $F_1:\Db(\cat{A})\to\cat{D}$ and
$F_2:\Db(\cat{A})\to\cat{D}$ be exact functors such that
\begin{itemize}
\item[(i)] there exists an isomorphism of functors
$f:F_2\rest{\cat{C}}\isomor F_1\rest{\cat{C}}$;
\item[(ii)] $\Hom_{\cat{D}}(F_1(A),F_1(B)[j])=0$,
for any $A,B\in\ob(\cat{A})$ and any $j<0$;
\item[(iii)] $F_1$ has a left adjoint $F_1^*$.
\end{itemize}
Then there exists an isomorphism of functors $g:F_2\isomor F_1$
extending $f$.\end{prop}

\begin{proof} We denote by $N$ the homological dimension of $\cat{A}$.

For any $i\in\ZZ$, let $f_i:=f(P_i):F_2(P_i)\isomor F_1(P_i)$. Given
$A\in\ob(\cat{A})$, we want to construct an isomorphism
$f_A:F_2(A)\isomor F_1(A)$. According to Lemma \ref{lem:resol}, let
\begin{eqnarray}\label{eqn:3}
\cdots\mor P_{i_j}^{\oplus k_j}\mor[d_j]P_{i_{j-1}}^{\oplus
k_{j-1}}\mor[d_{j-1}]\cdots\mor[d_1]P_{i_0}^{\oplus
k_0}\mor[d_0]A\mor0
\end{eqnarray}
be a resolution of $A$. Fix $m>N$ and consider the bounded complex
\[
R_m:=\{P_{i_m}^{\oplus k_m}\mor[d_m]P_{i_{m-1}}^{\oplus
k_{m-1}}\mor[d_{m-1}]\cdots\mor[d_1]P_{i_0}^{\oplus k_0}\}.
\]
Due to Remark \ref{rmk:sp}, a (unique up to isomorphism) convolution
of $R_m$ is $(d_0,0):P_{i_0}^{\oplus k_0}\to A\oplus K_m[m]$.

Due to Remark \ref{rmk:convfun}, for $i\in\{1,2\}$, the complex
\[
F_i(R_m):=\{F_i(P_{i_m}^{\oplus k_m})
\mor[F_i(d_m)]F_i(P_{i_{m-1}}^{\oplus k_{m-1}})\mor[F_i(d_{m-1})]
\cdots\mor[F_i(d_1)]F_i(P_{i_0}^{\oplus k_0})\}
\]
admits a convolution $(F_i(d_0),0):F_i(P_{i_0}^{\oplus k_0})\to
F_i(A\oplus K_m[m])$. Lemma \ref{lem:conv1} and conditions (i) and
(ii) ensure that such a convolution is unique up to
isomorphism. Moreover, again by (i) and (ii),
$\Hom_{\cat{D}}(F_2(P_{i_k}),F_1(A)[r])\iso
\Hom_{\cat{D}}(F_2(P_{i_l}),F_1(P_{i_j})[r])=0$, for any
$i_j,i_l,i_k\in\{i_0,\ldots,i_m\}$ and $r<0$. Hence we can apply
Lemma \ref{lem:conv2} getting a unique isomorphism
$f_A:F_2(A)\isomor F_1(A)$ making the following diagram commutative:
\[
\xymatrix{F_2(P_{i_m}^{\oplus k_m})\ar[r]^-{F_2(d_m)}
\ar[d]^{f_{i_m}^{\oplus k_m}}&F_2(P_{i_{m-1}}^{\oplus k_{m-1}})
\ar[r]^-{F_2(d_{m-1})}\ar[d]^{f_{i_{m-1}}^{\oplus k_{m-1}}}&
\cdots\ar[r]^-{F_2(d_1)}&F_2(P_{i_0}^{\oplus k_0})\ar[r]^-{F_2(d_0)}
\ar[d]^{f_{i_0}^{\oplus k_0}}&F_2(A)\ar[d]^{f_A}\\
F_1(P_{i_m}^{\oplus k_m})\ar[r]^-{F_1(d_m)}&F_1(P_{i_{m-1}}^{\oplus
k_{m-1}})\ar[r]^-{F_1(d_{m-1})}&\cdots\ar[r]^-{F_1(d_1)}&
F_1(P_{i_0}^{\oplus k_0})\ar[r]^-{F_1(d_0)}&F_1(A).}
\]

By Lemma \ref{lem:conv2}, the definition of $f_A$ does not depend on
the choice of $m$. In other words, if we choose a different $m'>N$
and we truncate \eqref{eqn:3} in position $m'$, the bounded
complexes $F_i(R_{m'})$ give rise to the same isomorphism $f_A$.

To show that the definition of $f_A$ does not depend on the choice
of the resolution \eqref{eqn:3}, consider another resolution
\begin{eqnarray}\label{eqn:reso2}
\cdots\mor P_{i'_j}^{\oplus k'_j}\mor[d'_j]P_{i'_{j-1}}^{\oplus
k'_{j-1}}\mor[d'_{j-1}]\cdots\mor[d'_1]P_{i'_0}^{\oplus
k'_0}\mor[d'_0] A\mor0.
\end{eqnarray}
Suppose that there exists a third resolution
\begin{eqnarray}\label{eqn:raff}
\cdots\mor P_{i''_j}^{\oplus k''_j}\mor[d''_j]P_{i''_{j-1}}^{\oplus
k''_{j-1}}\mor[d''_{j-1}]\cdots\mor[d''_1]P_{i''_0}^{\oplus
k''_0}\mor[d''_0] A\mor0
\end{eqnarray}
and morphisms $s_j:P_{i''_j}^{\oplus k''_j}\to P_{i_j}^{\oplus
k_j}$ and $t_j:P_{i''_j}^{\oplus k''_j}\to P_{i'_j}^{\oplus
k'_j}$, for any $j\geq 0$, fitting into the following commutative
diagram:
\[
\xymatrix{\cdots\ar[r]^{d'_{j+1}}&P_{i'_j}^{\oplus
k'_j}\ar[r]^{d'_j} &P_{i'_{j-1}}^{\oplus k'_{j-1}}\ar[r]^{d'_{j-1}}
&\cdots\ar[r]^{d'_1} & P_{i'_0}^{\oplus k'_0}\ar[dr]^{d'_0}&
\\
\cdots\ar[r]^{d''_{j+1}}&P_{i''_j}^{\oplus
k''_j}\ar[u]^{t_j}\ar[d]^{s_j}\ar[r]^{d''_j} & P_{i''_{j-1}}^{\oplus
k''_{j-1}}\ar[u]^{t_{j-1}}\ar[d]^{s_{j-1}}\ar[r]^{d''_{j-1}}
&\cdots\ar[r]^{d''_1} & P_{i''_0}^{\oplus
k''_0}\ar[d]^{s_{0}}\ar[u]^{t_{0}}\ar[r]^{d''_0}&
A.\\
\cdots\ar[r]^{d_{j+1}}&P_{i_j}^{\oplus k_j}\ar[r]^{d_j}
&P_{i_{j-1}}^{\oplus k_{j-1}}\ar[r]^{d_{j-1}} &\cdots\ar[r]^{d_1} &
P_{i_0}^{\oplus k_0}\ar[ur]_{d_0}& }
\]
Define the bounded complexes
\begin{eqnarray}\label{eqn:4}
\begin{array}{c}
R''_m:=\{P_{i''_m}^{\oplus k''_m}\mor[d''_m]P_{i''_{m-1}}^{\oplus
k''_{m-1}}\mor[d''_{m-1}]\cdots\mor[d''_1]P_{i''_0}^{\oplus k''_0}\}\\
F_i(R''_{m}):=\{F_i(P_{i''_m}^{\oplus k''_m})\mor[F_i(d''_m)]
F_i(P_{i''_{m-1}}^{\oplus k''_{m-1}})\mor[F_i(d''_{m-1})]
\cdots\mor[F_i(d''_1)]F_i(P_{i''_0}^{\oplus k''_0})\}.
\end{array}
\end{eqnarray}
Let $f''_A:F_2(A)\isomor F_1(A)$ be the isomorphism constructed
using \eqref{eqn:4}. Due to Remark \ref{rmk:sp}, these complexes and
their convolutions give rise to the diagram
\[
\xymatrix{F_2(P_{i''_0}^{\oplus k''_0})\ar[rrr]^{F_2(d''_0)}
\ar[dr]^{F_2(s_0)}\ar[ddd]^{f_{i''_0}^{\oplus k''_0}}& & &F_2(A)
\ar[dl]^{\id}\ar[ddd]^{f''_A}\\
&F_2(P_{i_0}^{\oplus k_0})\ar[r]^{F_2(d_0)}\ar[d]^{f_{i_0}^{\oplus
k_0}}
&F_2(A)\ar[d]^{f_A}&\\
&F_1(P_{i_0}^{\oplus k_0})\ar[r]^{F_1(d_0)}
&F_1(A)\ar@{}[ur]|{\bigstar}&\\
F_1(P_{i''_0}^{\oplus
k''_0})\ar[rrr]^{F_1(d''_0)}\ar[ur]_{F_1(s_0)}& & &
F_1(A)\ar[ul]_{\id}}
\]
where all squares but $\bigstar$ are commutative. Due to hypotheses
(i), (ii) and Lemma \ref{lem:conv2} there exists a unique morphism
$F_2(A)\to F_1(A)$ making the following diagram commutative:
\[
\xymatrix{F_2(P_{i''_0}^{\oplus
k''_0})\ar[d]_{F_1(s_0)\comp f_{i''_0}^{\oplus
k''_0}}\ar[rr]^{F_2(d''_0)}& &
F_2(A)\ar[d]\\
F_1(P_{i_0}^{\oplus k_0})\ar[rr]^{F_1(d_0)}& & F_1(A).}
\]
Since $F_1(s_0)\comp f_{i''_0}^{\oplus k''_0}=f_{i_0}^{\oplus
k_0}\comp F_2(s_0)$, both $f_A$ and $f''_A$ have this property and
then they coincide. Similarly one can prove that the morphism
$f''_A$ is equal to the morphism $f'_A$ constructed by means of
\eqref{eqn:reso2}.

To construct \eqref{eqn:raff}, we proceed as follows. First take
$i''_0\ll 0$ such that there exist a surjective morphism
$d''_0:P_{i''_0}^{\oplus k''_0}\epi A$, for some $k''_0\in\NN$, and
two morphisms $s_0$ and $t_0$ as required. Suppose now that
$P_{i''_j}$, $k''_j$, $d''_j$, $s_j$ and $t_j$ are defined. Take
$i''_{j+1}\ll 0$, $k''_{j+1}\in\NN$ and
$d''_{j+1}:P_{i''_{j+1}}^{\oplus k''_{j+1}}\to P_{i''_j}^{\oplus
k''_j}$ such that
\begin{itemize}
\item[(a.1)] $\ker(d''_j)=\im(d''_{j+1})$;
\item[(b.1)] the morphism $s_j\comp d''_{j+1}$ factorizes through
$d_{j+1}$;
\item[(c.1)] the morphism $t_j\comp d''_{j+1}$ factorizes through
$d'_{j+1}$.
\end{itemize}
Observe that this is always possible because
$\im(s_j\rest{\ker(d''_j)})\subset\im(d_{j+1})$ and for $n\ll0$ the
natural map $\Hom_{\cat{A}}(P_n,P^{\oplus
k_{j+1}}_{i_{j+1}})\to\Hom_{\cat{A}}(P_n,\im(d_{j+1}))$ is
surjective (the same holds true for $d'_{j+1}$ and $t_j$).

To prove the functoriality, let $A,B\in\ob(\cat{A})$ and let
$\varphi:A\to B$ be a morphism. Consider a resolution
\begin{eqnarray}\label{eqn:5}
\cdots\mor P_{l_j}^{\oplus h_j}\mor[e_j]P_{l_{j-1}}^{\oplus
h_{j-1}}\mor[e_{j-1}]\cdots\mor[e_1]P_{l_0}^{\oplus h_0}\mor[e_0]
B\mor0.
\end{eqnarray}
Reasoning as before, we can find a resolution
\begin{eqnarray}\label{eqn:6}
\cdots\mor P_{i_j}^{\oplus k_j}\mor[d_j]P_{i_{j-1}}^{\oplus
k_{j-1}}\mor[d_{j-1}]\cdots\mor[d_1]P_{i_0}^{\oplus
k_0}\mor[d_0]A\mor0
\end{eqnarray}
and morphisms $g_j:P_{i_{j}}^{\oplus k_j}\to P_{l_j}^{\oplus h_j}$ defining a
morphism of complexes compatible with $\varphi$. Fix $m>N$ and take the bounded
complexes
\[
\begin{array}{c}
R_m:=\{P_{i_m}^{\oplus k_m}\mor[d_m]P_{i_{m-1}}^{\oplus
k_{m-1}}\mor[d_{m-1}]\cdots\mor[d_1]P_{i_0}^{\oplus k_0}\}\\
T_m:=\{P_{l_m}^{\oplus h_m}\mor[e_m]P_{l_{m-1}}^{\oplus
h_{m-1}}\mor[e_{m-1}]\cdots\mor[e_1]P_{l_0}^{\oplus h_0}\}\\
F_i(R_{m}):=\{F_i(P_{i_m}^{\oplus k_m})\mor[F_i(d_m)]
F_i(P_{i_{m-1}}^{\oplus k_{m-1}})\mor[F_i(d_{m-1})]\cdots
\mor[F_i(d_1)]F_i(P_{i_0}^{\oplus k_0})\}\\
F_i(T_{m}):=\{F_i(P_{l_m}^{\oplus h_m})\mor[F_i(e_m)]
F_i(P_{l_{m-1}}^{\oplus h_{m-1}})\mor[F_i(e_{m-1})]\cdots
\mor[F_i(e_1)]F_i(P_{l_0}^{\oplus h_0})\}.
\end{array}
\]
We can now consider the diagram
\[
\xymatrix{F_2(P_{i_0}^{\oplus k_0})\ar[rrr]^{F_2(d_0)}
\ar[dr]^{f_{i_0}^{\oplus k_0}}
\ar[ddd]^{F_2(g_0)}& & & F_2(A)\ar[dl]^{f_A}\ar[ddd]^{F_2(\varphi)}\\
& F_1(P_{i_0}^{\oplus k_0})\ar[r]^{F_1(d_0)}\ar[d]^{F_1(g_0)}
&F_1(A)
\ar[d]^{F_1(\varphi)}&\\
&F_1(P_{l_0}^{\oplus h_0})\ar[r]^{F_1(e_0)}
&F_1(B)\ar@{}[ur]|{\bigstar}&\\
F_2(P_{l_0}^{\oplus h_0})\ar[rrr]^{F_2(e_0)}
\ar[ur]_{f_{l_0}^{\oplus h_0}}& & &F_2(B)\ar[ul]_{f_B}}
\]
where all squares but $\bigstar$ are commutative. Applying (i), (ii)
and Lemma \ref{lem:conv2} we see that there is a unique morphism
$F_2(A)\to F_1(B)$ completing the following diagram to a commutative
square
\[
\xymatrix{F_2(P_{i_0}^{\oplus k_0})\ar[d]_{F_1(g_0)\comp
f_{i_0}^{\oplus k_0}}\ar[rr]^{F_2(d_0)}& &
F_2(A)\ar[d]\\
F_1(P_{l_0}^{\oplus h_0})\ar[rr]^{F_1(e_0)}& & F_1(B).}
\]
Since $F_1(g_0)\comp f_{i_0}^{\oplus k_0}=f_{l_0}^{\oplus
h_0}\comp F_2(g_0)$, both $F_1(\varphi)\comp f_A$ and $f_B\comp
F_2(\varphi)$ have this property. Thus $F_1(\varphi)\comp
f_A=f_B\comp F_2(\varphi)$.

If $A\in\ob(\cat{A})$, we can clearly put $f_{A[n]}:=f_A[n]$ for
every integer $n$. Moreover, for any $A,B\in\ob(\cat{A})$, the
morphisms $f_A$ and $f_B$ just constructed commute with any
$g\in\Hom_{\Db(\cat{A})}(A,B[j])$ (see \cite[Sect.\ 2.16.4]{Or1} for
the proof).

The rest of the proof follows the strategy in \cite[Sect.\
2.16.5]{Or1} and it proceeds by induction on the length of the
segment in which the cohomologies of the objects are concentrated.
In particular, let $A$ be an object in $\Db(\cat{A})$ and suppose,
without loss of generality, that $H^p(A)=0$ if $p\not\in[a,0]$ and
$a<0$. Consider a morphism $v:P_i^{\oplus k}\to A$ such that
\begin{itemize}
\item[(a.2)] the natural morphism $u:P_i^{\oplus k}\to H^0(A)$
induced by $v$ is surjective; \item[(b.2)]
$\Hom_{\cat{A}}(H^0(F_1^*\comp F_1(A)),P_i)=0$.
\end{itemize}
Take a distinguished triangle
\[
Z[-1]\mor P_i^{\oplus k}\mor[v]A\mor Z
\]
and observe that $H^p(Z)=0$ if $p\not\in[a,-1]$. Hence, by
induction hypothesis, we have an isomorphism $f_Z:F_2(Z)\isomor
F_1(Z)$ and the following commutative diagram
\[
\xymatrix{F_1(Z)[-1]\ar[r]\ar[d]^{f^{-1}_Z[-1]}& F_1(P_i^{\oplus k})
\ar[r]\ar[d]^{(f_i^{\oplus k})^{-1}}&F_1(A)\ar[r]&F_1(Z)\ar[d]^{f^{-1}_Z}\\
F_2(Z)[-1]\ar[r]& F_2(P_i^{\oplus k})\ar[r]&F_2(A)\ar[r]&F_2(Z).}
\]

By \cite[Lemma 1.4]{Or1}, to complete the previous diagram with a
unique isomorphism $f_A:F_2(A)\isomor F_1(A)$, we need to show that
$$\Hom_{\cat{D}}(F_1(A),F_2(P_i))=0.$$ To prove this we can suppose $A\in
\ob(\cat{A})$ because the cohomologies of $A$ are concentrated in
degrees less or equal to zero. Let $w=\max\{n\in\ZZ:H^n(F^*_1\comp
F_1(A))\neq 0\}$. Obviously, there exists a natural non-zero
morphism $F^*_1\comp F_1(A)\to H^w(F^*_1\comp F_1(A))[-w]$. Hence
\[
0\neq\Hom_{\Db(\cat{A})}(F^*_1\comp F_1(A),H^w(F^*_1\comp
F_1(A))[-w])\iso\Hom_{\cat{D}}(F_1(A),F_1(H^w(F^*_1\comp F_1(A)))[-w])
\]
and $w\leq 0$ because of (ii). In particular $H^j(F^*_1\comp
F_1(A))=0$ if $j\not\in[-b,0]$, for some positive integer $b$.
Therefore, due to (b.2) and (ii),
\begin{multline*}
\Hom_{\cat{D}}(F_1(A),F_2(P_i))
\iso\Hom_{\cat{D}}(F_1(A),F_1(P_i))\\
\iso\Hom_{\Db(\cat{A})}(F_1^*\comp F_1(A),P_i)
\iso\Hom_{\cat{A}}(H^0(F_1^*\comp F_1(A)),P_i)=0.
\end{multline*}

To prove that $f_A$ is well-defined and functorial, one has to
repeat line by line the proof in Sections 2.16.6 and 2.16.7 of
\cite{Or1} using (ii) instead of the hypothesis that $F_1$ and
$F_2$ are fully-faithful. We leave this to the reader.\end{proof}

\begin{cor}\label{cor:ex} Let $\cat{A}$ and $\cat{B}$ be abelian categories
such that $\cat{A}$ has finite dimensional $\Hom$'s, it is of
finite homological dimension and it has an ample sequence. If
$F:\Db(\cat{A})\to\Db(\cat{B})$ is an exact functor with a left
adjoint and such that $F(\cat{A})\subseteq\cat{B}$, then
$G:=F|_{\cat{A}}:\cat{A}\to\cat{B}$ is exact and $\Db(G)\iso
F$.\end{cor}

\begin{proof} The exactness of $G$ is trivial. Since $\Db(G)|_{\cat{A}}\iso
F|_{\cat{A}}$, we can apply Proposition \ref{prop:extending} getting
the desired conclusion.\end{proof}

\section{Proof of Theorem \ref{thm:main}}

We divide up our argument in several steps.

\subsection{Resolution of the diagonal}\label{subsec:diag} Denoting by
$d:X\mono X\times X$ the diagonal morphism, $\ko_{\Delta}:=d_*\ko_X$ can be
regarded as an $(\alpha^{-1}\boxtimes\alpha)$-twisted coherent sheaf on
$X\times X$ in a natural way, since $d^*(\alpha^{-1}\boxtimes\alpha)=1$. It is
easy to see that $\ko_\Delta\in\coh(X\times X,\alpha^{-1}\boxtimes\alpha)$
admits a resolution
\begin{equation}\label{eqn:diag}
\cdots\mor A_i\boxtimes B_i\mor[\delta_i]A_{i-1}\boxtimes
B_{i-1}\mor[\delta_{i-1}]\cdots\mor[\delta_1]A_0\boxtimes
B_0\mor[\delta_0]\ko_\Delta\mor0,
\end{equation}
where $A_j\in\coh(X,\alpha^{-1})$ and $B_j\in\coh(X,\alpha)$ are
locally free for any $j\in\NN$. Indeed, if $L$ is an ample line
bundle on $X$, $L\boxtimes L$ is ample on $X\times X$. Hence, given
a locally free sheaf $E\in\coh(X,\alpha)$, Lemma \ref{lem:ample}
proves that $\{(E\dual\boxtimes E)\otimes(L\boxtimes L)^{\otimes
k}\}_{k\in\ZZ}$ is an ample sequence in $\coh(X\times
X,\alpha^{-1}\boxtimes\alpha)$. As $(E\dual\boxtimes
E)\otimes(L\boxtimes L)^{\otimes k}\iso(E\dual\otimes L^{\otimes
k})\boxtimes(E\otimes L^{\otimes k})$, we conclude by Lemma
\ref{lem:resol}.

\subsection{Some bounded complexes}\label{subsec:bc} Since $F$ is a bounded
functor by Proposition \ref{prop:bound}, we can assume without loss of
generality that $H^i(F(\kf))=0$ for any $\kf\in\coh(X,\alpha)$ and any
$i\notin[-M,0]$ for some $M\in\NN$. Then we fix once and for all a resolution
of $\ko_\Delta$ as in \eqref{eqn:diag}, and for every integer
$m>\dim(X)+\dim(Y)+M$ we define the following complexes
\begin{gather*}
C_m:=\{A_m\boxtimes B_m\mor[\delta_m]\cdots\mor[\delta_1]
A_0\boxtimes B_0\}\\
\tilde{C}_m:=\{A_m\boxtimes F(B_m)\mor[\tilde{\delta}_m]\cdots
\mor[\tilde{\delta}_1]A_0\boxtimes F(B_0)\}
\end{gather*}
in $\Db(X\times X,\alpha^{-1}\boxtimes\alpha)$ and $\Db(X\times Y,
\alpha^{-1}\boxtimes\beta)$ respectively, where $\tilde{\delta}_i$ denotes the
image of $\delta_i$ through the map
\begin{multline*}
\Hom_{\Db(X\times X,\alpha^{-1}\boxtimes\alpha)} (A_i\boxtimes
B_i,A_{i-1}\boxtimes B_{i-1})\iso \Hom_{\Db(X,\alpha^{-1})}(A_i,A_{i-1})\otimes
\Hom_{\Db(X,\alpha)}(B_i,B_{i-1})\\
\mor[\id\otimes F]\Hom_{\Db(X,\alpha^{-1})}(A_i,A_{i-1})\otimes
\Hom_{\Db(Y,\beta)}(F(B_i),F(B_{i-1}))\\
\iso\Hom_{\Db(X\times Y,\alpha^{-1}\boxtimes\beta)} (A_i\boxtimes
F(B_i),A_{i-1}\boxtimes F(B_{i-1})).
\end{multline*}
Setting $\kk_m:=\ker(\delta_m)\in\coh(X\times
X,\alpha^{-1}\boxtimes\alpha)$ and proceeding as in Remark
\ref{rmk:sp}, we see that, if $m>2\dim(X)$,
$\cp{C_m}\iso\ko_{\Delta}\oplus\kk_m[m]$ and $C_m$ has a (unique up
to isomorphism by Lemma \ref{lem:conv1}) right convolution
$(\delta_0,0):A_0\boxtimes B_0\to\ko_{\Delta}\oplus\kk_m[m]$.
Observe that the assumption on $F$ implies that also $\tilde{C}_m$
satisfies the hypothesis of Lemma \ref{lem:conv1}, hence it has a
unique up to isomorphism right convolution $\tilde{\delta}'_{0,m}:
A_0\boxtimes F(B_0)\to\kg_m$.

We proceed now as in \cite[Lemma 6.1]{Ka}. We denote by $\cat{K}_m$
the full subcategory of $\coh(X,\alpha)$ with objects the locally
free sheaves $E$ such that $H^i(X,E\otimes A_j)=0$ for $i>0$ and
$0\le j\le m+\dim(X)$. Observe that, for any locally free
$E'\in\coh(X,\alpha)$ and any ample line bundle $L\in\coh(X)$,
$E'\otimes L^{\otimes k}\in\cat{K}_m$, when $k\gg 0$.

As $\R^i{p_2}_*(A_j\boxtimes B_j\otimes p_1^*\kf)\iso
H^i(X,A_j\otimes\kf)\otimes B_j$ for $\kf\in\coh(X,\alpha)$ and
$i,j\in\NN$ (where $p_l:X\times X\to X$ is the projection onto the
$l^{\rm th}$ factor),
\[\R{p_2}_*(A_j\boxtimes B_j\otimes p_1^*E)\iso
{p_2}_*(A_j\boxtimes B_j\otimes p_1^*E)\iso H^0(X,A_j\otimes E)\otimes B_j\] if
$E\in\cat{K}_m$ and $0\le j\le m+\dim(X)$. It follows that the exact functor
$\R{p_2}_*(-\otimes p_1^*E)$ maps $C_m$ to a complex
\[C_{m,E}=\{H^0(X,A_m\otimes E)\otimes B_m\mor[\delta_{m,E}]\cdots
\mor[\delta_{1,E}]H^0(X,A_0\otimes E)\otimes B_0\}\] in $\Db(X,\alpha)$, which
has a (unique up to isomorphism) right convolution
\[(\delta_{0,E},0):H^0(X,A_0\otimes E)\otimes B_0\to E\oplus\kk_{m,E}[m],\]
where $\kk_{m,E}:=\ker(\delta_{m,E})$. If $p:X\times Y\to Y$ and
$q:X\times Y\to X$ are the natural projections, a similar argument
shows that the exact functor $\R p_*(-\otimes q^*E)$ maps
$\tilde{C}_m$ to a complex
\[\tilde{C}_{m,E}=\{H^0(X,A_m\otimes E)\otimes F(B_m)\mor[\tilde{\delta}_{m,E}]
\cdots\mor[\tilde{\delta}_{1,E}]H^0(X,A_0\otimes E)\otimes F(B_0)\}\] in
$\Db(Y,\beta)$, which has a unique up to isomorphism right convolution
\begin{equation}\label{eqn:E1}
\tilde{\delta}'_{0,m,E}:H^0(X,A_0\otimes E)\otimes F(B_0)\to\R p_*(\kg_m\otimes
q^*E)=\FM{\kg_m}(E).
\end{equation}
On the other hand, $\tilde{C}_{m,E}$ can be identified with the image of
$C_{m,E}$ through $F$, so that a right convolution of $\tilde{C}_{m,E}$ is
given also by
\begin{equation}\label{eqn:E2}
(F(\delta_{0,E}),0):H^0(X,A_0\otimes E)\otimes F(B_0)\to F(E)\oplus
F(\kk_{m,E})[m].
\end{equation}
Therefore $\FM{\kg_m}(E)\iso F(E)\oplus F(\kk_{m,E})[m]$, and so, in
particular, $H^i(\FM{\kg_m}(E))=0$ unless $i\in[-m-M,-m]\cup[-M,0]$.
Since this holds for every $E\in\cat{K}_m$, applying Lemma
\ref{lem:serrerel} we deduce that also $H^i(\kg_m)=0$ unless
$i\in[-m-M,-m]\cup[-M,0]$. This implies that
$\kg_m\iso\ke_m\oplus\kf_m$ with $H^i(\ke_m)=0$ unless $i\in[-M,0]$
and $H^i(\kf_m)=0$ unless $i\in[-m-M,-m]$.

\subsection{Uniqueness of the kernel} We are going to show that a
kernel of $F$ (if it exists) is necessarily isomorphic to $\ke_m$ for $m\gg0$.
Indeed, assume that $\ke\in\Db(X\times Y, \alpha^{-1}\boxtimes\beta)$ is such
that $F\iso\FM{\ke}$. A standard computation shows that
\[\ke':=p_{13}^*\ko_{\Delta}\lotimes p_{24}^*\ke\in
\Db(X\times X\times X\times Y,
\alpha\boxtimes\alpha^{-1}\boxtimes\alpha^{-1}\boxtimes\beta)\]
(where $p_{ij}$ is the obvious projection from $X\times X\times
X\times Y$ and $\ko_{\Delta}$ is now considered to be
$(\alpha\boxtimes\alpha^{-1})$-twisted) defines a functor of
Fourier-Mukai type
\[\FM{\ke'}:\Db(X\times X,\alpha^{-1}\boxtimes\alpha)\mor
\Db(X\times Y,\alpha^{-1}\boxtimes\beta)\] such that
$\FM{\ke'}(\kf)\iso\R(q_{13})_*(q_{12}^*\kf\lotimes q_{23}^*\ke)$ for every
$\kf\in \Db(X\times X,\alpha^{-1}\boxtimes\alpha)$ (here, again, $q_{ij}$
denotes the obvious projection from $X\times X\times Y$). It follows easily
that $\FM{\ke'}(\ko_{\Delta})\iso\ke$ and
\[\FM{\ke'}(\ka\lboxtimes\kb)\iso\ka\lboxtimes\FM{\ke}(\kb)\iso
\ka\lboxtimes F(\kb)\] for $\ka\in\Db(X,\alpha^{-1})$ and
$\kb\in\Db(X,\alpha)$. In particular, we see that $\FM{\ke'}$ maps
the complex $C_m$ to $\tilde{C}_m$, hence if $m>2\dim(X)$ a
convolution of the latter complex is given by
$\FM{\ke'}(\ko_{\Delta}\oplus\kk_m[m])\iso
\ke\oplus\FM{\ke'}(\kk_m)[m]$. Therefore
$\ke\oplus\FM{\ke'}(\kk_m)[m] \iso\kg_m\iso\ke_m\oplus\kf_m$, and we
can conclude that $\ke\iso\ke_m$ (and
$\FM{\ke'}(\kk_m)[m]\iso\kf_m$) provided $m\gg0$ (more precisely, it
is enough that $\Hom_{\Db(X\times
Y,\alpha^{-1}\boxtimes\beta)}(\ke,\kf_m)=0$ and $\Hom_{\Db(X\times
Y,\alpha^{-1}\boxtimes\beta)}(\ke_m,\FM{\ke'}(\kk_m)[m])=0$, which
is certainly true for large $m$ by definition of $\ke_m$ and $\kf_m$
and because $\FM{\ke'}$ is bounded).

\subsection{Isomorphism of functors on a subcategory} Now we fix an integer
$m>\dim(X)+\dim(Y)+M$ and we will prove that $\ke:=\ke_m$ is
really a kernel of $F$. To simplify the notation we will suppress
the subscript $m$ also from $\kg_m$, $\kf_m$, $\tilde{C}_{m,E}$,
$\tilde{\delta}'_{0,m,E}$ and $\cat{K}_m$. As a first step, we
will show that $\FM{\ke}\rest{\cat{K}}$ and $F\rest{\cat{K}}$ are
isomorphic as functors from $\cat{K}$ to $\Db(Y,\beta)$. To see
this we use the argument in \cite[Lemma 6.2]{Ka}. In fact for
every $E\in\cat{K}$ by \eqref{eqn:E1} and \eqref{eqn:E2} the
complex $\tilde{C}_E$ has two right convolutions, namely
\[\tilde{\delta}'_{0,E}=(\tilde{\delta}_{0,E},0):H^0(X,A_0\otimes E)\otimes
F(B_0)\to\FM{\kg}(E)\iso\FM{\ke}(E)\oplus\FM{\kf}(E)\] and
$(F(\delta_{0,E}),0)$. Due to Lemma \ref{lem:conv2} this implies
that there exists a unique isomorphism
$\varphi(E):\FM{\ke}(E)\isomor F(E)$ such that
$F(\delta_{0,E})=\varphi(E)\comp\tilde{\delta}_{0,E}$. In order to
see that this isomorphism is functorial, just notice that for every
morphism $\gamma:E\to E'$ of $\cat{K}$, again by Lemma
\ref{lem:conv2}, there is a unique morphism $\FM{\ke}(E)\to F(E')$
such that the diagram
\[\xymatrix{H^0(X,A_0\otimes E)\otimes F(B_0) \ar[rr]^-{\tilde{\delta}_{0,E}}
\ar[d]_{H^0(\id\otimes\gamma)\otimes\id} & & \FM{\ke}(E) \ar[d] \\
H^0(X,A_0\otimes E')\otimes F(B_0) \ar[rr]^-{F(\delta_{0,E'})} & &
F(E')}\] commutes. Since both $F(\gamma)\comp\varphi(E)$ and
$\varphi(E')\comp\FM{\ke}(\gamma)$ satisfy this property, they
must be equal.

\subsection{Extending the isomorphism} Now we choose an $\alpha$-twisted
locally free sheaf $E$ and a very ample line bundle $L$ on $X$
(defining an embedding $X\mono\PP^N$) and we denote by $\cat{C}$ the
full subcategory of $\coh(X,\alpha)$ with objects $\{E\otimes
L^{\otimes k}\}_{k\in\ZZ}$. Now, by Lemma \ref{lem:ample} this set
of objects is an ample sequence in $\coh(X,\alpha)$, hence by
Proposition \ref{prop:extending} in order to prove that
$F\iso\FM{\ke}$ it is enough to show that
$F\rest{\cat{C}}\iso\FM{\ke}\rest{\cat{C}}$. To this purpose, we
proceed as in \cite[Lemma 6.4]{Ka} and we define isomorphisms
$\varphi_k:F(E\otimes L^{\otimes k})\isomor\FM{\ke}(E\otimes
L^{\otimes k})$ (for $k\in\ZZ$) such that
\begin{equation}\label{eqn:nat}
\FM{\ke}(\gamma)\comp\varphi_{k_1}=\varphi_{k_2}\comp F(\gamma)
\end{equation}
for every morphism $\gamma:E\otimes L^{\otimes k_1}\to E\otimes
L^{\otimes k_2}$ of $\cat{C}$ and for every $k_1,k_2\in\ZZ$. By
definition of ample sequence there exists $k_0\in\ZZ$ such that
$E\otimes L^{\otimes k}\in\cat{K}$ for $k\ge k_0$. Then, setting
$\varphi_k:=\varphi(E\otimes L^{\otimes k})^{-1}$ for $k\ge k_0$,
the equation \eqref{eqn:nat} is satisfied for $k_1,k_2\ge k_0$. Now
we proceed by descending induction: assuming $\varphi_k$ is defined
for $k>n$ and \eqref{eqn:nat} is satisfied for $k_1,k_2>n$, we
define $\varphi_n$ as follows. As in \eqref{eqn:Beilinson},
Beilinson's resolution gives an exact sequence in $\coh(X)$
\[
0\mor\ko_X\mor[\rho_{N+1}]L\otimes V_N\mor[\rho_N]\cdots\mor[\rho_2]
L^{\otimes N}\otimes V_1\mor[\rho_1]L^{\otimes N+1}\otimes V_0\mor0
\]
(where each $V_i$ is a finite dimensional vector space), hence,
setting $\rho^{(n)}_i:=\id_{E\otimes L^{\otimes n}}\otimes\rho_i$,
the complex
\[E\otimes L^{\otimes n+1}\otimes V_N\mor[\rho^{(n)}_N]\cdots\mor[\rho^{(n)}_1]
E\otimes L^{\otimes n+N+1}\otimes V_0\] in $\Db(X,\alpha)$ has a unique up to
isomorphism left convolution $\rho^{(n)}_{N+1}:E\otimes L^{\otimes n}\to
E\otimes L^{\otimes n+1}\otimes V_N$. The inductive hypothesis implies that
\[\xymatrix{F(E\otimes L^{\otimes n+1})\otimes V_N \ar[rr]^-{F(\rho^{(n)}_N)}
\ar[d]^{\varphi_{n+1}\otimes\id} & & \cdots \ar[rr]^-{F(\rho^{(n)}_1)} & &
F(E\otimes L^{\otimes n+N+1})\otimes V_0 \ar[d]^{\varphi_{n+N+1}\otimes\id} \\
\FM{\ke}(E\otimes L^{\otimes n+1})\otimes V_N \ar[rr]^-{\FM{\ke}(\rho^{(n)}_N)}
& & \cdots \ar[rr]^-{\FM{\ke}(\rho^{(n)}_1)} & & \FM{\ke}(E\otimes L^{\otimes
n+N+1})\otimes V_0}\] is an isomorphism of complexes in $\Db(Y,\beta)$ which
satisfies the assumptions of Lemma \ref{lem:conv2}, hence there is a unique
isomorphism $\varphi_n$ such that the diagram
\[\xymatrix{F(E\otimes L^{\otimes n}) \ar[rr]^-{F(\rho^{(n)}_{N+1})}
\ar[d]^{\varphi_n} & & F(E\otimes L^{\otimes n+1})\otimes V_N
\ar[d]^{\varphi_{n+1}\otimes\id} \\
\FM{\ke}(E\otimes L^{\otimes n}) \ar[rr]^-{\FM{\ke}(\rho^{(n)}_{N+1})} & &
\FM{\ke}(E\otimes L^{\otimes n+1})\otimes V_N}\] commutes. Moreover, for every
morphism $\gamma:E\otimes L^{\otimes k_1}\to E\otimes L^{\otimes k_2}$ of
$\cat{C}$ and for every $k_1,k_2\ge n$
\[\xymatrix{F(E\otimes L^{\otimes k_1+1})\otimes V_N
\ar[rr]^-{F(\rho^{(k_1)}_N)} \ar[d]^{\tilde{\gamma}_N\otimes\id} & & \cdots
\ar[rr]^-{F(\rho^{(k_1)}_1)} & & F(E\otimes L^{\otimes k_1+N+1})\otimes V_0
\ar[d]^{\tilde{\gamma}_0\otimes\id} \\
\FM{\ke}(E\otimes L^{\otimes k_2+1})\otimes V_N
\ar[rr]^-{\FM{\ke}(\rho^{(k_2)}_N)} & & \cdots
\ar[rr]^-{\FM{\ke}(\rho^{(k_2)}_1)} & & \FM{\ke}(E\otimes L^{\otimes
k_2+N+1})\otimes V_0}\] (where $\tilde{\gamma}_i:=\varphi_{k_2+N+1-i}\comp
F(\gamma\otimes\id_{L^{\otimes N+1-i}})=\FM{\ke}(\gamma\otimes\id_{L^{\otimes
N+1-i}})\comp\varphi_{k_1+N+1-i}$) is a morphism of complexes in $\Db(Y,\beta)$
which again satisfies the assumptions of Lemma \ref{lem:conv2}. Therefore,
there is a unique morphism $F(E\otimes L^{\otimes k_1})\to\FM{\ke}(E\otimes
L^{\otimes k_2})$ such that the diagram
\[\xymatrix{F(E\otimes L^{\otimes k_1}) \ar[rr]^-{F(\rho^{(k_1)}_{N+1})}
\ar[d] & & F(E\otimes L^{\otimes k_1+1})\otimes V_N
\ar[d]^{\tilde{\gamma}_N\otimes\id} \\
\FM{\ke}(E\otimes L^{\otimes k_2}) \ar[rr]^-{\FM{\ke}(\rho^{(k_2)}_{N+1})} & &
\FM{\ke}(E\otimes L^{\otimes k_2+1})\otimes V_N}\] commutes, and, since both
$\FM{\ke}(\gamma)\comp\varphi_{k_1}$ and $\varphi_{k_2}\comp F(\gamma)$ satisfy
this property, we conclude that \eqref{eqn:nat} holds.

\begin{remark}\label{Rmk:kaw} Theorem 1.1 in
\cite{Ka} concerns fully faithful functors. This requirement is
essential in Kawamata's proof only in \cite[Lemma 6.5]{Ka} (which
depends on \cite[Prop.\ 2.16]{Or1}). Kawamata's argument can now be
reconsidered using Proposition \ref{prop:extending} instead of
\cite[Lemma 6.5]{Ka}. Hence we immediately get the following
generalization of Kawamata's result. Let $X$ and $Y$ be normal
projective varieties with only quotient singularities and let $\kx$
and $\ky$ be the smooth stacks naturally associated to them. Let
$F:\Db(\coh(\kx))\to\Db(\coh(\ky))$ be an exact functor with a left
adjoint and such that, for any $\kf,\kg\in\coh(\kx)$,
$$\Hom_{\Db(\coh(\ky))}(F(\kf),F(\kg)[j])=0$$ if $j<0$. Then there exists a
unique up to isomorphism $\ke\in\Db(\coh(\kx\times\ky))$ and an isomorphism of
functors $F\cong\FM{\ke}$.

Observe moreover that, in Kawamata's proof, the results in
\cite[Sect.\ 3]{Ka} can be replaced by our shorter argument in
Section \ref{subsec:diag}.\end{remark}

\section{Exact functors between the abelian categories of twisted
sheaves}\label{sec:exfun}

Theorem \ref{thm:main} can be used to classify exact functors from
$\coh(X,\alpha)$ to $\coh(Y,\beta)$.

\begin{prop}\label{prop:exab}
Let $(X,\alpha)$ and $(Y,\beta)$ be twisted varieties. If $\ke$ is
in $\coh(X\times Y,\alpha^{-1}\boxtimes\beta)$, then the additive
functor
\[\aFM{\ke}:=p_*(\ke\otimes q^*(-)):\coh(X,\alpha)\to\coh(Y,\beta)\] is exact
if and only if $\ke$ is flat over $X$ and $p\rest{\supp(\ke)}:\supp(\ke)\to Y$
is a finite morphism.

Moreover, for every exact functor $G:\coh(X,\alpha)\to\coh(Y,\beta)$ there
exists unique up to isomorphism $\ke\in\coh(X\times
Y,\alpha^{-1}\boxtimes\beta)$ (flat over $X$ and with $p\rest{\supp(\ke)}$
finite) such that $G\iso\aFM{\ke}$.
\end{prop}

\begin{proof}
Clearly $\ke$ is flat over $X$ if and only if the functor
$\ke\otimes q^*(-)$ is exact, and in this case $\aFM{\ke}$ is left
exact and $\R\aFM{\ke}\iso\R p_*(\ke\otimes q^*(-))$. Notice also
that $\ke\otimes q^*(-)$ is exact if $\aFM{\ke}$ is exact. Indeed,
given an injective morphism $\kf\mono\kg$ in $\coh(X,\alpha)$ and
setting $$\kk:=\ker(\ke\otimes q^*\kf\to\ke\otimes q^*\kg),$$ for
every $E\in\coh(X)$ locally free there is an exact sequence in
$\coh(Y,\beta)$
\[0\to p_*(\kk\otimes q^*E)\to
p_*(\ke\otimes q^*(\kf\otimes E))=\aFM{\ke}(\kf\otimes E)\to p_*(\ke\otimes
q^*(\kg\otimes E))=\aFM{\ke}(\kg\otimes E),\] from which we see that, if
$\aFM{\ke}$ is exact, $p_*(\kk\otimes q^*E)=0$; therefore $\kk=0$, and this
proves that $\ke\otimes q^*(-)$ is exact. It follows that, in order to conclude
the proof of the first statement, it is enough to show that
$p\rest{\supp(\ke)}$ is finite if and only if $\R^j p_*(\ke\otimes q^*\kf)=0$
for $j>0$ and for every $\kf\in\coh(X,\alpha)$. To this purpose, up to
replacing $\ke$ with $\ke\otimes p^*F$ for some locally free sheaf
$F\in\coh(Y,\beta^{-1})$, we can assume that $\beta$ is trivial (because
$\supp(\ke)=\supp(\ke\otimes p^*F)$ and $\R^j p_*(\ke\otimes p^*F\otimes
q^*\kf)\iso F\otimes\R^j p_*(\ke\otimes q^*\kf)$).

If $p\rest{\supp(\ke)}$ is finite, we can find a cover by open affine subsets
$\{V_i\}_{i\in I}$ of $Y$ and $U_i\subset X$ open affine such that
$\supp(\ke\rest{X\times V_i})\subset U_i\times V_i$ for every $i\in I$. Then,
denoting by $p_i:X\times V_i\to V_i$ and $p'_i:U_i\times V_i\to V_i$ the
projections, for $j>0$ and for every $\kf\in\coh(X,\alpha)$ we have
\[\R^jp_*(\ke\otimes q^*\kf)\rest{V_i}\iso
\R^j{p_i}_*((\ke\otimes q^*\kf)\rest{X\times V_i})\iso
\R^j{p'_i}_*((\ke\otimes q^*\kf)\rest{U_i\times V_i})=0\] because
$p'_i$ is an affine morphism, hence $\R^jp_*(\ke\otimes
q^*\kf)=0$. On the other hand, if $p\rest{\supp(\ke)}$ is not
finite, there exist a closed point $y\in Y$ and a closed
irreducible subset $X'\subseteq X$ such that $d:=\dim(X')>0$ and
$X'\subseteq\supp(\ke_y)$, where $\ke_y\in\coh(X,\alpha^{-1})$
corresponds to $\ke\rest{X\times\{y\}}$ under the natural
isomorphism $X\iso X\times\{y\}$. We claim that there exists
$\kf_0\in\coh(X,\alpha)$ such that $\supp(\kf_0)=X'$ and
$H^d(X,\ke_y\otimes\kf_0)\ne0$. For instance, denoting by $E$ a
locally free $\alpha$-twisted sheaf on $X$ and by $\ko_{X'}(1)$ a
very ample line bundle on $X'$ (regarded as a subscheme of $X$ with
the reduced induced structure), we can take
$\kf_0=E\otimes\ko_{X'}(-n)$ for $n\gg0$. Indeed, by definition of
the dualizing sheaf $\ds_{X'}$ (see \cite[p.\ 241]{Ha}), we have
\begin{multline*}
H^d(X,\ke_y\otimes E\otimes\ko_{X'}(-n))\dual\iso H^d(X',(\ke_y\otimes
E)\rest{X'}(-n))\dual \\
\iso\Hom_{X'}((\ke_y\otimes E)\rest{X'},\ds_{X'}(n))\iso
H^0(\sHom_{X'}((\ke_y\otimes E)\rest{X'},\ds_{X'})(n)),
\end{multline*}
and the last term is not $0$ for $n\gg0$, since
$\sHom_{X'}((\ke_y\otimes E)\rest{X'},\ds_{X'})\ne0$ due to the
fact that $\supp((\ke_y\otimes E)\rest{X'})=X'$ and
$\ds_{X'}\iso\omega_{X'}$ on the non-empty open subset where $X'$ is smooth.

Then let $V\subset Y$ be an open affine subset containing $y$ and denote by
$q':X\times V\to X$ the projection. Applying the right exact functor
$q'_*(-)\otimes\kf_0$ to the natural surjective morphism $\ke\rest{X\times
V}\epi\ke\rest{X\times}\{y\}$, we get a surjective morphism in $\qcoh(X)$
\[\varphi:q'_*(\ke\rest{X\times V})\otimes\kf_0\epi
q'_*(\ke\rest{X\times\{y\}})\otimes\kf_0\iso\ke_y\otimes\kf_0.\] As
$\supp(\ker(\varphi))\subseteq X'$, we have
$H^{d+1}(X,\ker(\varphi))=0$, hence the assumption
$H^d(X,\ke_y\otimes\kf_0)\ne0$ implies that
\[0\ne H^d(X,q'_*(\ke\rest{X\times V})\otimes\kf_0)\iso
H^d(X\times V,(\ke\otimes q^*\kf_0)\rest{X\times V}),\] and this proves that
$\R^dp_*(\ke\otimes q^*\kf_0)\ne0$.

Assume now that $G:\coh(X,\alpha)\to\coh(Y,\beta)$ is an exact
functor. By Theorem \ref{thm:main} there exists (unique up to
isomorphism) $\ke\in\Db(X\times Y,\alpha^{-1}\boxtimes\beta)$ such
that $\Db(G)\iso\FM{\ke}$, and $\ke\in\coh(X\times
Y,\alpha^{-1}\boxtimes\beta)$ by Lemma \ref{lem:serrerel}. From the
fact that $\FM{\ke}(\coh(X,\alpha))\subseteq\coh(Y,\beta)$ it is
easy to deduce that
$$G\iso\FM{\ke}\rest{\coh(X,\alpha)}\iso\aFM{\ke}.$$ The uniqueness
of $\ke$ follows from Corollary \ref{cor:ex}.
\end{proof}

\begin{remark}
The above result implies that there are no non-zero exact functors
from $\coh(X,\alpha)$ to $\coh(Y,\beta)$ if $\dim(X)>\dim(Y)$: to
prove this, just note that if $0\ne\ke\in\coh(X\times
Y,\alpha^{-1}\boxtimes\beta)$ is flat over $X$ then
$\dim(\supp(\ke))\ge\dim(X)$, and that
$\dim(\supp(\ke))\le\dim(Y)$ if $p\rest{\supp(\ke)}$ is finite.
\end{remark}

It was proved by Gabriel in \cite{G} that if $X$ and $Y$ are
noetherian schemes then there exists an exact equivalence
$\qcoh(X)\iso\qcoh(Y)$ if and only if $X$ and $Y$ are isomorphic.
For smooth projective varieties, a short proof (relying on Orlov's
result) of an analogous statement involving coherent sheaves was
given in \cite{Hu1}. Following this last approach and using
Proposition \ref{prop:exab} we prove a Gabriel-type result for
twisted varieties.

\begin{cor} Let $(X,\alpha)$ and $(Y,\beta)$ be twisted varieties.
Then the following three conditions are equivalent:
\begin{itemize}
\item[(i)] there is an exact equivalence between $\qcoh(X,\alpha)$ and
$\qcoh(Y,\beta)$;
\item[(ii)] there is an exact equivalence between $\coh(X,\alpha)$ and
$\coh(Y,\beta)$;
\item[(iii)] there exists an isomorphism $f:X\isomor Y$ such that
$f^*(\beta)=\alpha$.
\end{itemize}\end{cor}

\begin{proof} The implications $(\mathrm{iii})\Rightarrow(\mathrm{i})$ and
$(\mathrm{iii})\Rightarrow(\mathrm{ii})$ are trivial. Suppose that
an exact equivalence $G:\qcoh(X,\alpha)\isomor\qcoh(Y,\beta)$ is
assigned and consider the equivalence
$\D(G):\D(\qcoh(X,\alpha))\isomor\D(\qcoh(Y,\beta))$ induced by $G$.
Due to \cite[Thm.\ 18]{P} and \cite[Lemma 2.1.4, Prop.\ 2.1.8]{C},
$\D(G)$ restricts to an equivalence
$F:\Db(X,\alpha)\isomor\Db(Y,\beta)$ which yields an exact
equivalence $G':\coh(X,\alpha)\isomor\coh(Y,\beta)$. This proves
that $(\mathrm{i})$ implies $(\mathrm{ii})$.

The proof of the implication
$(\mathrm{ii})\Rightarrow(\mathrm{iii})$ proceeds now as in
\cite[Cor.\ 5.22, Cor.\ 5.23]{Hu1}. First of all, recall that given
an abelian category $\cat{A}$, $A\in\ob(\cat{A})$ is \emph{minimal}
if any non-trivial surjective morphism $A\to B$ in $\cat{A}$ is an
isomorphism. Notice that an equivalence $F:\cat{A}\isomor\cat{B}$
sends minimal objects to minimal objects.

It is easy to see that the set of minimal objects of
$\coh(X,\alpha)$ consists of all skyscraper sheaves $\ko_x$, where
$x$ is a closed point of $X$. Let
$G:\coh(X,\alpha)\isomor\coh(Y,\beta)$ be an exact equivalence. By
Proposition \ref{prop:exab}, $G\iso\aFM{\ke}$, for some
$\ke\in\coh(X\times Y,\alpha^{-1}\boxtimes\beta)$.

Since $G$ maps skyscraper sheaves to skyscraper sheaves,
$\ke\rest{\{x\}\times Y}$ is isomorphic to a skyscraper sheaf and we
naturally get a morphism $f:X\to Y$ and $L\in\Pic(X)$ such that
$$G\iso L\otimes f_*(-).$$ The morphism $f$ is actually an isomorphism
since $G$ is an equivalence (\cite[Cor.\ 5.23]{Hu1}) and, by
definition, $f^*(\beta)=\alpha$.\end{proof}

This proves Conjecture 1.3.17 in \cite{C} for quasi-coherent sheaves
on smooth projective varieties.

\begin{remark} Suppose that $X$ and $Y$ are smooth
separated schemes of finite type over a field $\K$ and that
$\alpha\in\Br(X)$ while $\beta\in\Br(Y)$. In \cite{P} it was proved
that if there exists an exact equivalence
$\coh(X,\alpha)\iso\coh(Y,\beta)$, then there is an isomorphism
$f:X\isomor Y$. On the other hand the approach in \cite{P} does not
allow to conclude that $f^*(\beta)=\alpha$.\end{remark}

\medskip

{\small\noindent{\bf Acknowledgements.} This paper was partially
written while the second named author was visiting the
Max-Planck-Institut f\"{u}r Mathematik (Bonn). It is a great
pleasure to thank the staff for the exciting mathematical atmosphere
at the Institute.}

\end{document}